\documentclass[aps,preprint,groupedaddress,longbibliography]{revtex4}

\usepackage{graphicx,dcolumn,bm}
\usepackage{amssymb,amstext,amsmath}
\usepackage{graphicx}
\usepackage{dcolumn}
\usepackage{bm}
\usepackage{amssymb}
\usepackage{multirow}
\usepackage{bigstrut}
\usepackage{makecell,rotating}
\usepackage{mathrsfs}
\usepackage{booktabs}
\usepackage{threeparttable}
\usepackage{multirow}
\usepackage{subfigure}
\usepackage{epsfig}
\usepackage{threeparttable}
\usepackage{chngpage}
\usepackage{float}
\usepackage{amstext}
\usepackage{amsmath}
\usepackage{setspace}
\usepackage{amsfonts}

\newcommand{\x}{\mathbf{x}}
\newcommand{\A}{\mathbf{A}}
\newcommand{\B}{\mathbf{B}}

\newcommand{\W}{\mathbf{W}}

\newcommand{\zero}{\mathbf{0}}

\renewcommand{\u}{\mathbf{u}}

\newcommand{\Lc}{\mathcal{L}}
\newcommand{\Rc}{\mathcal{R}}

\newcommand{\e}{\textrm{e}}

\newcommand{\ie}{\textit{i.e.}}

\begin{document}

\title{The optimal trajectory to control complex networks}
\author{Aming Li$^{1,}$\footnote{amingli2011@gmail.com}, 
Long Wang$^{2,}$\footnote{longwang@pku.edu.cn}, 
and Frank Schweitzer$^{1,}$\footnote{fschweitzer@ethz.ch}}

\affiliation{
\begin{enumerate}
  \item Chair of Systems Design, Department of Management, 
Technology and Economics, ETH Z\"urich, Weinbergstrasse 56/58, Z\"urich CH-8092, Switzerland
  \item Center for Systems and Control, College of Engineering, Peking University, Beijing 100871, China
\end{enumerate}
}
\date{\today}

\begin{abstract}
Controllability, a basic property of various networked systems, has gained profound theoretical applications in complex social, technological, biological, and brain networks.
Yet, little attention has been given to the control trajectory (route), along which a controllable system can be controlled from any initial to any final state, hampering the implementation of practical control.
Here we systematically uncover the fundamental relations between control trajectory and several other key factors, such as the control distance between initial and final states ($\delta$), number of driver nodes, and the control time.
The length ($\Lc$) and maximum distance to the initial state ($\Rc$) are employed to quantify the locality and globality of control trajectories. 
We analyze how the scaling behavior of the averaged $\Lc$ and $\Rc$ changes with increasing $\delta$ for different initial states. 
After showing the scaling behavior for each trajectory, we also provide the distributions of $\Lc$ and $\Rc$.
Further attention is given to the control time $t_{f}$ and its influence on $\Lc$ and $\Rc$. 
Our results provide comprehensive insights in understanding control trajectories for complex networks, and pave the way to achieve practical control in various real systems.
\end{abstract}

\maketitle

\newpage
\section{Introduction}
As a powerful framework, complex networks have been widely employed to understand various complex systems,
where nodes indicate system's components and links capture interactions between them \cite{SocialNet94book,Havlin2004book,Barabasi04,barrat2008dynamical,Rajapakse2011pnas,Coyte2015Sci}.
Controllability---a basic property detecting whether a system can be controlled from external inputs 
\cite{Liu2011,Wang2013,Gao2014,Posfai2014NJP,LiXiang2014,Liu2016Rev,Cornelius2013NatCommun,Chen2014,Chen2017,Li2017,Gu08Con,ConBrain2016Simu,Yan2017Nature,Bassett18NP}, helps to 
uncover the principles of, for example, the interactions of neural circuits of cognitive function in brain networks \cite{Gu08Con}, or even predicting neuron function in the nematode \textit{Caenorhabditis elegans} \cite{Yan2017Nature}.
Indeed, a system is said to be controllable, if it can be driven from arbitrary initial state to arbitrary final state within finite time under appropriate control inputs  \cite{Kalman63,Xie2002tac,Chen2017}.
However, the reported principles of control cannot tell how systems behave under control inputs, namely, no information can be obtained on the evolution of a system's state in order to reach the desired state by just testing the system's controllability.
Although some results emerge on control cost (energy)   
\cite{Rajapakse2011pnas,Yan2012PRL,Yan2015a,Lai2016,Klickstein2017,Bassett18NP}, the practical control trajectory (route) from the initial to final state along which the system must traverse all transient states, is far from understood, which strongly inhibits the practical applications.

Here we systematically explore control trajectories for controlling complex networks, revealing the fundamental relations between practical trajectories and control distance, number of driver nodes, and the control time.
Our findings clarify the fundamental behavior of practical control trajectories when we control complex networks, impulsing the real applications of the network control theory.

\section{Dynamics on Complex Networks}
The dynamics of a complex network with external inputs can be described mathematically as 
\begin{equation}
\label{nonlineardy}
\dot{\mathcal{X}}(t) = \boldsymbol{f}(t, \mathcal{X} (t), \mathcal{U}(t), \mathcal{P}),
\end{equation}
where $\mathcal{X}_i(t)$ is the state of node $i$ at time $t$, like the level of neural activity of brain region $i$ in a brain network \cite{Gu08Con,Bassett18NP,Yan2017Nature},
or the concentration of metabolite $i$ in a metabolic network \cite{Almaas2004}.
The vector $\mathcal{X}(t)$ collects the state of all the $N$ nodes, \ie, $\mathcal{X}(t) = (\mathcal{X}_1(t),\mathcal{X}_2(t),\cdots,\mathcal{X}_N(t))^\textrm{T} \in\mathbb{R}^N$, represents the system state at time $t$.
$\boldsymbol{f}(*) = (f_{1}(*),f_{2}(*),\cdots,f_{N}(*))^\textrm{T}$ denotes interaction dynamics among nodes.
$\mathcal{U}(t)\in\mathbb{R}^p$ captures the input signals acting directly on $p$ ($\leq N$) nodes (namely, driver nodes \cite{Liu2011}).
$\mathcal{P}$ is the set of the system's parameters, which reflects the exact intensity that nodes interact with each other. 

\begin{figure}
\centering
\includegraphics[width=1\textwidth]{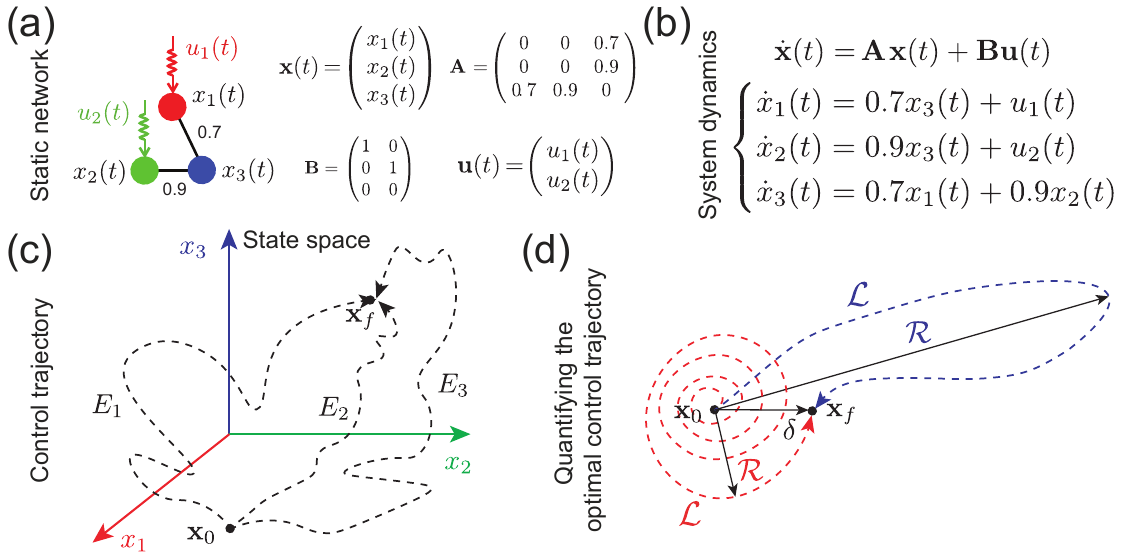}
\caption{
Networks and the related dynamics.
In (a), we show a network with 3 nodes for clarity. 
We employ $x_i(t)$ to represent the state of node $i$ at time $t$, and hence $\x(t)$ represents the state of the whole network.
The corresponding adjacency matrix is given aside the network.
The matrix $\B$ gives the mapping between inputs $\u(t)$ and the driver nodes, which receive inputs directly as shown in red and green nodes. 
The dynamics described in equation (\ref{lineardy}) is presented in (b), where it shows how the state of each node evolves under the control inputs given in a.
In the system's state space of $\x(t)$ plotted in (c), we denote the initial and final states of the network in (a) as $\x_0$ and $\x_f$.
With appropriate control inputs $\u(t)$, we can drive the system's state from $\x_0$ to $\x_f$. 
For different  $\u(t)$, there are different trajectories, among which we show $3$ different ones, and the corresponding control cost $E_1, E_2$, and $E_3$ are given aside. 
Among all the possible control trajectories starting from $\x_0$ to $\x_f$, here we focus the optimal one along which the control cost is minimal. 
(d), Two typical variables to quantify the optimal control trajectory, one is the length ($\Lc$) showing how long the trajectory wanders totally until reaching the final state, and another is the radius ($\Rc$) telling the longest distance the trajectory reaches from the initial state.
$\delta$ is the direct control distance between the initial and final state.
}
\label{fig_1_cartoon}
\end{figure}

Due to the lack of empirical information about the 
exact nonlinearity of $\boldsymbol{f}(*)$ and the related set of precise parameters $\mathcal{P}$, equation (\ref{nonlineardy}) is normally linearized to pursue analytical insights \cite{Liu2011,Gao16Nature,Yan2017Nature,Gu2017Traj,Bassett18NP}.
By assuming that the fixed point of the network is $\mathcal{X}^*$ without additional inputs, \ie, $\boldsymbol{f}(t, \mathcal{X}^*, \mathcal{U}^*)=0$, we linearize (\ref{nonlineardy}) by employing $\x(t) = \mathcal{X}(t) - \mathcal{X}^*$ and $\u(t) = \mathcal{U}(t)-\mathcal{U}^*$, arriving at the following dynamics
\begin{equation}
\label{lineardy}
\dot{\x }(t) = \A \x(t) + \B\u(t),
\end{equation}
in the time interval $[t_0, t_f]$ (see Fig.~\ref{fig_1_cartoon}a).
$\A= \frac{\partial \boldsymbol{f}(*)}{\partial  \mathcal{X}} \Big |_{\mathcal{X}^*, \mathcal{U}^*}$ 
corresponds to the adjacency matrix of the network (see Fig.~\ref{fig_1_cartoon}b and c), whose entry $a_{ij}$ represents, for example, the number of white matter streamlines linking from regions $j$ to $i$ in the brain network \cite{Gu2017Traj,Bassett18NP}. 
$\B = \frac{\partial  \boldsymbol{f}(*)}{\partial  \mathcal{U} }\Big |_{\mathcal{X}^*, \mathcal{U}^*}$ gives the constant mapping between inputs and driver nodes of the network (see Fig.~\ref{fig_1_cartoon}a).

\section{Variables to quantify the control trajectory}

To quantify the control trajectory when we control complex networks, we adopt two variables.
One is the length 
\begin{equation}
\label{Ldefmain}
\Lc = \int_{t_0}^{t_f}\| \dot{\x }(t)\| \mathrm{d}t
= \int_{t_0}^{t_f} \sqrt{ \dot{\x }^{\textrm{T}}(t) \dot{\x }(t) }\mathrm{d}t
\end{equation}
telling how long the control trajectory wanders in the controllable space. 
Indeed, the length of control trajectory is widely used to quantify the locality of control trajectories for complex networks \cite{Sun2013prl,Li2017}, and it is also applied to analyze brain networks \cite{Gu2017Traj}.
It is discovered that $\Lc$ can be extremely large when the control distance $\delta = \|\x_f - \x_0 \|$ approaches $0$ \cite{Sun2013prl,Li2017}.
This implies that the optimal trajectory is probably nonlocal where in some dimensions the state components of the trajectory pass through highly extreme values (see Fig.~\ref{fig_1_cartoon}d).
Nevertheless, when $\Lc$ is large, it does not necessarily mean that the optimal trajectory is nonlocal.
Indeed, when the trajectory circuits around the initial state before arriving at the final one, $\Lc$ can still be large but the system state does not wander far from the initial state (see Fig.~\ref{fig_1_cartoon}d).
It means that the optimal trajectory cannot be solely reflected by the magnitude of $\Lc$. 
Here we propose the radius of the control trajectory  
\begin{eqnarray}
\label{DefDmain}
\Rc &=& \max_{t_0 \leq t \leq t_f} \| \x(t) - \x_0 \| 
=
\max_{t_0 \leq t \leq t_f} \sqrt{ \sum_{i=1}^N \left(x_i(t) - x_i(t_0) \right)^2}
\end{eqnarray}
to quantify the maximum distance that the control trajectory deviates from the initial state among all of the system's intermediate states.
Here $\Rc$ can serve as a signal to dictate the existence of extreme values of state components. 
Indeed, if there are some extremely large values of $x_i(t)$, then $\Rc$ will be large as well, and if the control trajectory is direct from the initial to the final state, then we have $\Rc \approx \delta$.

\section{The optimal control trajectory}

For the dynamics given in equation~(\ref{lineardy}), we obtain that, starting from $\x_0$ at time $t_0$, the control trajectory at the time $t$ $(<t_f)$ is
\begin{eqnarray}
\x (t)  =  \textrm{e}^{\A (t-t_0)}\x _0 + \int^t_{t_0} \textrm{e}^{\A (t-\tau)}\B \u(\tau) \textrm{d} \tau,
\end{eqnarray}
with the external input $\u(\tau)$.
To drive the network to reach the final state $\x_f$ at time $t_f$, however, we can choose an enormous number of different inputs (Fig.~\ref{fig_1_cartoon}d), which in turn generate different control trajectories with different control costs.
Indeed, the input control cost is defined as $E = \int^{t_f}_{t_0} \u(t)^{\text{T}} \u(t) \text{d} t$ \cite{OptimalBooLewis}, which reaches its minimum with the optimal control input 
 $$\u(t)=\B^{\textrm{T}}\e^{\A ^\textrm{T}(t_{f}-t)}\W^{-1}[t_0, t_f] \textbf{\textrm{d}} $$
where 
$\textbf{\textrm{d}} = \x_f - \e^{\A (t_f-t_0)}\x_0 $ is the difference between the desired final state $\x_f$ and the natural final state that the system evolves without external inputs, and 
$\W = \int^{t_f}_{t_0} \e^{\A (t_f-\tau)}\B \B^{\textrm{T}}\e^{\A ^\textrm{T}(t_{f}-\tau)}  \textrm{d} \tau.$
Here, for given initial and final states, we focus on the optimal control trajectory determined by the optimal control inputs, along which the control cost is minimum.

\section{How initial states and control distances affect the averaged length (radius) of control trajectories}

When the network is steered from $\x_0$ to $\x_f$ in practice, it is of great interest how the direct control distance $\delta = \|\x_f - \x_0 \|$ affects the way from $\x_0$ to $\x_f$.
The length of the optimal control trajectory is
\begin{eqnarray}
\Lc =
 \int_{t_0}^{t_f}  \sqrt{ 
\| \x_0\|^2f(\bar{\x}_0,\bar{\x}_0) + 2\| \x_0\| \| \x_f\| f(\bar{\x}_0,\bar{\x}_f) + \| \x_f\|^2f(\bar{\x}_f,\bar{\x}_f)
}  \mathrm{d}t,
\end{eqnarray}
where $\bar{\x}_0$ and $\bar{\x}_f$ is the unit vector along the direction of $\x_0$ and $\x_f$ separately, and the function $f(*)$ is given in the Ref.~\cite{SM}.
The final state $\x_f = \x_0 + \delta \bar{\x}$ when $\bar{\x}$ is the unit vector along the direction of $\x_f - \x_0$.
This suggests that the behavior of control trajectories is determined by the relation between the initial state and the control distance.
Here we first focus on the overall behavior of the averaged length ($\Lc$) of control trajectories under the same direct control distance as a function of $\delta$.

\begin{figure}
\centering
\includegraphics[width=1\textwidth]{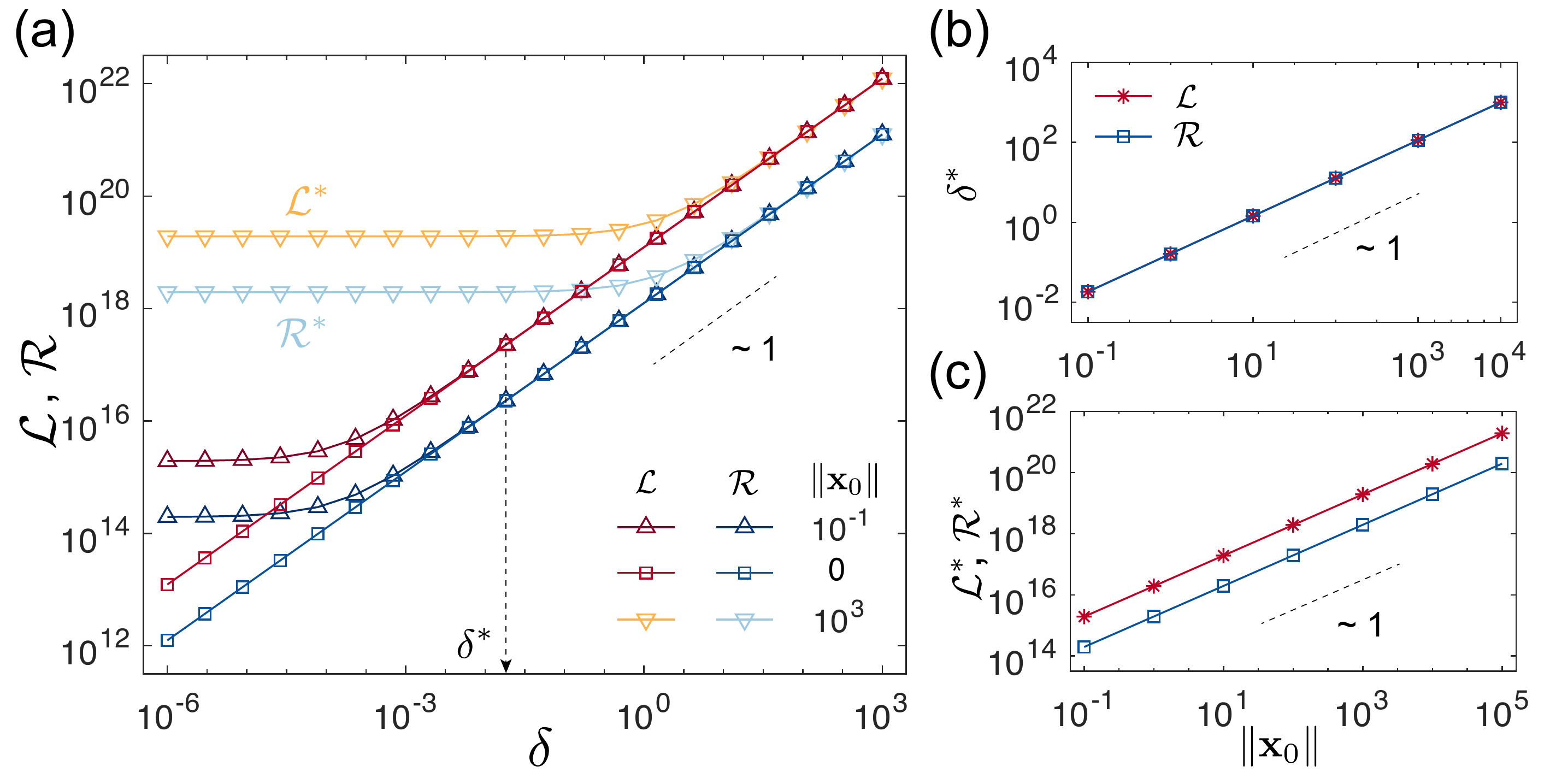}
\caption{
Scaling behavior of the averaged length and radius of control trajectories under different control distances.
(a), For each value of control distance $\delta$, different locations of $\x_0$ ($10^{-1},~0,~\text{and}~10^3$) are chosen to calculate the length $\Lc$ and radius $\Rc$ of the optimal control trajectories along which the control energy is minimum.
We choose $100$ final states $\x_f$ randomly on the sphere centered on $\x_0$ with the distance $\delta = \| \x_f - \x_0 \|$, over which the averaged $\Lc$ and $\Rc$ are obtained.
The scaling behavior of $\Lc$ for $\| \x_0 \| = 0$ is $\Lc \sim \delta$.
For $\| \x_0 \| \neq 0$, the scaling behavior of $\Lc$ depends on the competition between the magnitude of $\| \x_0 \|$ and $\delta$, where $\Lc$ first keeps as a constant $\Lc^*$ determined by the nonzero initial state $\x_0$ when $\delta < \delta^*$, and then is dominated by $\delta$ when $\delta > \delta^*$.
(b), We further find that for the critical value of the control distance $\delta^*$, at which the scaling behavior of $\Lc$ and $\Rc$ alters, linearly increases with the magnitude of the initial state $\| \x_0 \|$.
(c), As to the constants $\Lc^*$ and $\Rc^*$, they increase with $\| \x_0 \|$ as well with the scaling behavior $ \Lc^*\sim \| \x_0 \|$.
All of the above results have been approximated by analytical derivations \cite{SM}.
$\delta^*$ is calculated as the minimal $\delta$ which makes the distance between $\Lc, \Rc$ for $\x_0 \neq \zero$ and $\Lc, \Rc$ for $\x_0 = \zero$ smaller than $10^{-2}$.
Here $N=7$ with the average degree $4$, and the number of driver nodes is $1$.
We choose $40$ points along the control trajectory to numerically approximate $\Lc$ and $\Rc$.
For other values of the related parameters, please see Figs.~S1 and~S2.
}
\label{fig_L_delta}
\end{figure}

When $\x_0 = \zero$ ($\|\x_f\| = \delta$), we know that 
$\Lc(\zero, l \x_f) 
=
l \Lc(\zero, \x_f)$.
That means, when a network is controlled from the origin, the averaged length of the control trajectory increases linearly with the control distance, \ie, $\Lc \sim \delta$ (see Fig.~\ref{fig_L_delta}a).

When $\x_0 \neq \zero$, from $\x_f = \x_0 + \delta \bar{\x}$, we find that:
(i) With the increase of $\delta$ (say, bigger than the critical value $\delta^*$), the effect of $\x_0$ can be neglected, leading to  $\Lc(\x_0, l\x_f)  \approx l \Lc(\zero, \x_f)$,
which follows the laws of the scenario for $\x_0 = \zero$.
That is to say, when the control distance is relatively long compared to the norm of the initial state, it will dominate the scaling behavior of the averaged length of control trajectories (see Fig.~\ref{fig_L_delta}a);
(ii) When the control distance is relatively short ($\delta<\delta^*$) with a nonzero initial state, we find that the averaged $\Lc$ can be approximated by the constant
\begin{eqnarray}
\label{Lstarx0}
\Lc^*= 
\|\x_0 \| \int_{t_0}^{t_f}  \sqrt{ 
f(\bar{\x}_0,\bar{\x}_0) + 2 f(\bar{\x}_0,\bar{\x}_f) + f(\bar{\x}_f,\bar{\x}_f)
}  \mathrm{d}t, 
\end{eqnarray}
This means that the averaged length of the control trajectory is dominated by $\| \x_0 \|$ as a constant when the control distance is short (see Fig.~\ref{fig_L_delta}a).

Equation (\ref{Lstarx0}) also tells us that the averaged constant increases linearly with the norm of the initial state, \ie, $\Lc^* \sim \| \x_0 \|$ since $\Lc^*(l\x_0,\x_f) = l\Lc^*(\x_0,\x_f)$ (see Fig.~\ref{fig_L_delta}c).

As to the critical value of $\delta$ at which the behavior of the averaged $\Lc$ will alter, we know that when $\x_0 = \zero$, $\Lc = k_1 \delta$, and when $\x_0 \neq \zero$, the corresponding constant is $k_2 \| \x_0 \|$, both $k_1$ and $k_2$ are constants.
Thus, at the critical control distance $\delta^*$,  we have $k_1 \delta^* = k_2 \| \x_0 \|$, meaning that the scaling behavior of $\delta^*$ follows $\delta^* \sim  \| \x_0 \|$.
This is also validated with numerical calculations (see Fig.~\ref{fig_L_delta}b).

Taken together, we find an universal linear scaling behavior of both the averaged length and averaged radius of the optimal control trajectory, namely, $\Lc \left(\Rc\right) \sim \delta$, $\delta^* \sim  \| \x_0 \|$, and $\Lc^* \left(\Rc^*\right) \sim \| \x_0 \|$.

\section{Scaling behavior of each control trajectory and its distribution}
\begin{figure}
\centering
\includegraphics[width=0.7\textwidth]{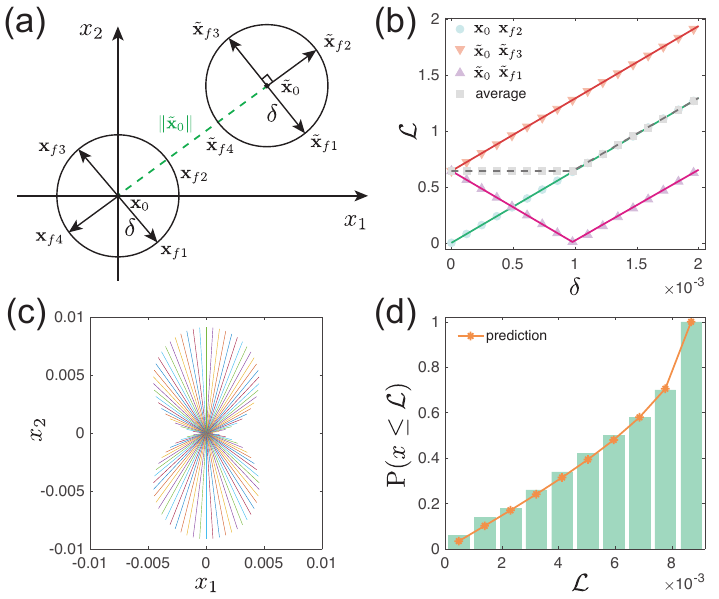}
\caption{
Scaling behavior and distribution of the length of control trajectories.
(a), Schematic presentation of the locations of initial and final states.
For a two dimensional system, we show two control scenarios with initial state at the origin ($\x_0 = \zero$) and away from the origin ($\x_0 = \tilde{\x}_0 \neq \zero$).
(b), Scaling behavior of three control trajectories as a function of the direct control distance $\delta$.
The green solid circles correspond to $\x_0 = \zero$ and $\x_f = \x_{f2}$ shown in (a), and the corresponding line generated from linear fitting, which shows $\Lc = a\delta$ with $R^2 = 1$.
For $\x_0 = \tilde{\x}_0$, the solid upward-pointing triangles represent $|Lc$ for $\x_f = \tilde{\x}_{f1}$, where the linear fitting gives $\Lc = -a\delta + b$ when $\delta<10^{-3}$, and $\Lc = a\delta - b$ for the rest.
As to the opposite direction ($\x_f = \tilde{\x}_{f3}$), results are presented in solid downward-pointing triangles, where we have $\Lc = a\delta + b$.
The averaged $\Lc$ over the cases for $\tilde{\x}_{f1}$ and $\tilde{\x}_{f1}$ is shown in grey solid square, which first keeps as a constant and then shares the same scaling law with that for $\x_0 = \zero$.
Here $a = 647.10, b = 0.64$, and the lines are generated from linear fitting of the corresponding dots with $R^2 > 0.999$.
Results for other control directions are given in Fig.~S4.
(c), Length of $100$ control trajectories for short control distance ($\delta=10^{-5}$) when $ \x_0  = \zero$.
Following each control direction $\x_{fi}$ selected uniformly ($i = 1,2,\cdots,100$), we plot the straight line with the length of the corresponding $\Lc$.
(d), The accumulated distribution of $\Lc$ shown in (c), where the solid line represents the analytical prediction from Eq.~(\ref{mainadfsi}).
Here the control time $t_f = 10^{-2}$, and the system is given in Eq.~(S12) \cite{SM}.
Other parameters are the same as those in Fig.~\ref{fig_L_delta}.
The results for other parameters, and the similar quantitative behavior of $\Rc$ are presented in Figs.~S4 and~S5.
Robust results for higher systems are given in Figs.~S6 to~S8.
 }
\label{fig_dist_LR}
\end{figure}

The averaged values of $\Lc$ and $\Rc$ provide statistical insights of control trajectories at the same control distance.
In the phase space, however, for two opposite final states ($\x_{f1}$ and $\x_{f3}$ in Fig.~\ref{fig_dist_LR}a), when their control distances to a given initial state ($\x_0$ in Fig.~\ref{fig_dist_LR}) are equal, they can correspond to totally different control objectives.
Indeed, for neural activity ($x_i(t_f)$) of the brain region $i$, the two final states $x_i(t_f)=1$ and $0$ have the same distance to the initial state $x_i(t_f)=0.5$, but $1$ and $0$ capture totally opposite states.
Thus, simply averaging over $\Lc$ or $\Rc$ for trajectories with the same $\delta$ may probably miss out the potential fundamental laws behind the practical control routes.
To better understand this, we first focus on each separate trajectory and then explore the statistical characteristics of all trajectories.   

Interestingly, we find that for nonzero initial state, $\Lc$ has the inverse scaling behavior for the opposite final states with the same small $\delta$.
For example, when $\x_0 = \tilde{\x}_0 \neq \zero$ (Fig.~\ref{fig_dist_LR}a), we randomly select a final state ($\x_f = \tilde{\x}_{f1}$) with direct distance $\delta$ to $\tilde{\x}_0$.
We find that the corresponding length of control trajectory first decreases with $\Lc = -a\delta + b$ and then increases linearly with $\Lc = a\delta - b$ (solid upward-pointing triangle in Fig.~\ref{fig_dist_LR}b).
As to the opposite direction ($\x_f = \tilde{\x}_{f3}$), we have $\Lc = a\delta + b$ (solid downward-pointing triangle in Fig.~\ref{fig_dist_LR}b).
When we average $\Lc$ over the final states with $\tilde{\x}_{f1}$ and $\tilde{\x}_{f3}$, we find that $\Lc$ first stays constant and then shares the same scaling law as for $\x_0 = \zero$ (grey solid square in Fig.~\ref{fig_dist_LR}b), which is in line with the results reported in Fig.~\ref{fig_L_delta}.
Thus the averaged $\Lc$ over different control trajectories with same control distance neutralizes the inverse scaling behavior for opposite final sates. 

For $\x_0 = \zero$ (Fig.~\ref{fig_dist_LR}a) and $\x_f = \x_{f2}$, we have $\Lc = a\delta$ (green solid circle in Fig.~\ref{fig_dist_LR}b), and $\Lc \sim \delta$ holds for any specific final state \cite{SM}.
In addition, as to any pair of opposite final states, the lengths of control trajectories are equal \cite{SM}.
Furthermore, for all the control trajectories at the same control distance (Fig.~\ref{fig_dist_LR}c), the cumulative distribution function of $\Lc$ is
\begin{eqnarray}
P (\Lc \leq x) = \frac{2}{\pi}  \arcsin \frac{x}{2r},
\label{mainadfsi}
\end{eqnarray}
where $2r$ is the maximum value of $\Lc$.
The above function can predict the numerical results very well (Fig.~\ref{fig_dist_LR}d). 

When $\x_0 \neq \zero$, the constant nonzero initial state determines the uniform distribution of $\Lc$ for small $\delta$, while for large $\delta$, $\Lc$ has the same distribution given by the above equation for both zero and nonzero initial states.
All the above results are applicable for the radius of control trajectories ($\Rc$), and other more results are given in the Ref.~\cite{SM}.

\section{How control time affects the scaling behavior of $\Lc$ and $\Rc$}
Under a given control distance, the control time ($t_f-t_0$) that control signals can harness to drive the system to the final state is quite important.
It affects not only the velocity of system state change but also the corresponding minimal control energy.
Here we seek to address how the control time affects the scaling behavior of $\Lc$ and $\Rc$.
According to equation (\ref{Ldefmain}), we have
\begin{eqnarray}
\label{Lbasicfort}
\Lc 
=
\int_{0}^{t_f} \sqrt{
\x_f ^\text{T} \W^{-1}[0, t_f] \textrm{e}^{\A(t_{f}-t)} \left( \W[0,t]\A  + \mathbf{I} \right)
\left( \A \W[0,t] + \mathbf{I} \right)
\textrm{e}^{\A ^\textrm{T}(t_{f}-t)}
\W^{-1}[0, t_f] \x_f 
 }\mathrm{d}t.
 \end{eqnarray}
To theoretically analyze the relation between $\Lc$ and the control time $t_f$, we divide it into three situations according to the number of driver nodes, \ie, one driver node, $p~(1 < p <N)$ driver nodes, and $N$ driver nodes \cite{SM}.
Note that, without loss of generality, here we set $\x_0 = \zero$ and $t_0 = 0$.
We numerically show the results as follows.

\begin{figure}
\centering
\includegraphics[width=1\textwidth]{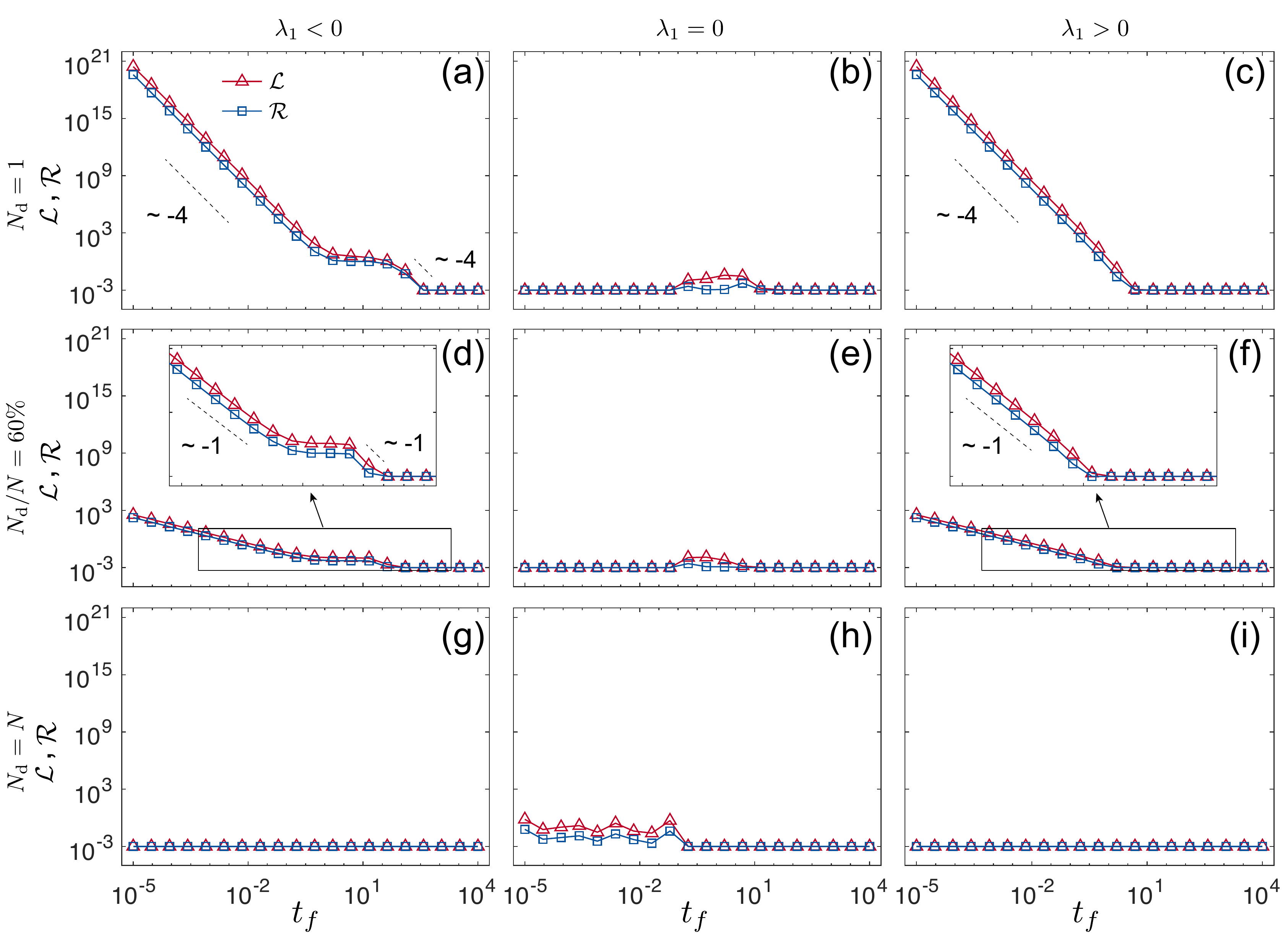}
\caption{
Scaling behavior of the length and radius of control trajectories under different control time.
For different numbers of driver nodes, the scaling behavior of $\Lc$ and $\Rc$ is determined by the largest eigenvalue ($\lambda_1$) of the adjacency matrix $\A$.
(a), With one driver node ($N_{\text{d}} = 1$), when $\A$ is negative definite ($\lambda_1<0$), we find that the scaling behavior of $\Lc$ and $\Rc$ decreases $t_f^{-4}$ for short control time $t_f$.
With the increase of $t_f$ $\Lc$ and $\Rc$ will first keep as a constant then decrease again as a power-law function of the middle level $t_f$ (from $10^1$ to $10^2$).
And eventually $\Lc$ and $\Rc$ keep as a constant for large control time.
(b), When $\A$ is negative semi-definite ($\lambda_1=0$), both $\Lc$ and $\Rc$ keep as a constant.
(c), When $\A$ is not negative definite ($\lambda_1>0$), $\Lc$ and $\Rc$ will first decrease with $t_f^{-4}$ for small $t_f$ and then keep as a constant as $t_f$ is large.
The increase of the driver nodes diminishes both $\Lc$ and $\Rc$, while maintaining the type of the scaling ((d) to (f)).
Indeed, as $N_{\text{d}}/N = 60 \%$, we find that the scaling behavior of $\Lc$ and $\Rc$ is $t_f^{-1}$ for short time when $\lambda_1 \neq 0$ (insets of panels d and f).
When we control all the nodes directly ($N_{\text{d}} = N$),  $\Lc$ and $\Rc$ keep as a constant, where control time cannot diminish the control trajectories ((g) to (i)).
All of the above results have been approximated by analytical derivations \cite{SM}.
Here $N=5$ with the average degree $3.5$, $\x_0 = \zero$ and
$\delta=10^{-3}$.
Other parameters are the same as those in Fig.~\ref{fig_L_delta}. 
For other values of the related parameters, please see Fig.~S9.
}
\label{fig4_N5_LR_t}
\end{figure}

For one driver node and short control time $t_f$, we find that $\Lc$ ($\Rc$) decreases with the power-law function of $t_f$ when the system is asymptotically stable (the maximum eigenvalue $\lambda_1$ of $\A$ is smaller than $0$) or unstable ($\lambda_1>0$) (Fig.~\ref{fig4_N5_LR_t}a and~\ref{fig4_N5_LR_t}c).
With the increase of $t_f$ for an asymptotically stable system, $\Lc$ ($\Rc$) will first keep as a constant and then decrease again with the same power-law function, and eventually keep as the constant $\alpha$ for big $t_f$ (Fig.~\ref{fig4_N5_LR_t}a).
We find that $\alpha \approx \delta$ and $\Lc \approx \Rc$ if $t_f$ is big, meaning that the control trajectory goes straight from the initial to the final state when the control time is long enough.

Interestingly, when the system is unstable, $\Lc$ ($\Rc$) keeps $\delta$ with the increase of $t_f$ (Fig.~\ref{fig4_N5_LR_t}c), while in this case we know that the minimum control energy is $E_{\min} \sim \text{e}^{- 2 \lambda_1 t_f}$ \cite{Yan2012PRL,Li2017ConEng}.
That is to say, for unstable systems, when more control \emph{time} is given, the corresponding optimal trajectory stays constant despite that the minimum energy needed to reach final state decreases exponentially. 
For the critical scenario where $\lambda_1 = 0$, we find that $\Lc$ ($\Rc$) equals  $\delta$ irrespective of how much control time is given (Fig.~\ref{fig4_N5_LR_t}b).
This means that, although the control time is short, an increasing control time can reduce the control energy dramatically \cite{Yan2012PRL,Li2017ConEng}.
But neither the length nor the radius of the control trajectory can be secured.

By adding more driver nodes, both $\Lc$ and $\Rc$ decrease, and the exponent of the scaling behavior of $\Lc$ and $\Rc$ will decrease as well (Fig.~\ref{fig4_N5_LR_t}d-f).
When we control all nodes directly, \ie, when the number of driver nodes is equal to the system size, both the length and radius of control trajectories keep constant for different scenarios of stability of the system and control time (Fig.~\ref{fig4_N5_LR_t}g-i).

\section{Conclusion and Discussions}

We statistically analyze the averaged length and the averaged radius of control trajectories with the same control distance.
We also provide the scaling behavior of these two quantities. 
We demonstrate that aggregating the length (radius) of trajectories over many evenly selected final states \emph{neutralizes} the embedded scaling behavior for each single trajectory.
For example, as $\x_0 = \zero$, the linear scaling behavior of $\Lc$ ($\Rc$) for every pair of opposite final states has the contrary sign for short control distance. Averaging this will make $\Lc$ ($\Rc$) a constant.
Thus the statistical results are not enough to fully understand the control trajectories.
Apart from uncovering the relations of the scaling for different final states equidistant to $\x_0$, we also analytically provide the distribution of $\Lc$ and $\Rc$.
In addition, $\Lc$ and $\Rc$ can be employed to classify different kinds of optimal control trajectory in terms of the locality and globality in various empirical systems. 

Another key factor to implement control under practical circumstances is the control time ($t_f$), \ie~the time needed to reach the final state.
We find that for short $t_f$, $\Lc$ ($\Rc$) is a power law function of $t_f$, meaning that, in this case $\Lc$ ($\Rc$) can be dramatically reduced if slightly more time is given.
When $t_f$ is big, $\Lc$ ($\Rc$) cannot be affected too much either by increasing the number of driver nodes or by changing the stability of the system.
This has consequences \textit{e.g.} for cognitive control, where the brain can quickly achieve some complex cognitive functions by altering the dynamics of neural systems with energetic inputs \cite{CogCon02,CogCon13,CogCon14}. 
Our findings suggest that the $\Lc$ ($\Rc$) of the optimal control trajectories in the phase space of neural activity can be largely conserved when more time is given to the brain to perform the cognitive control.

To pursue the analytical insights of optimal control trajectories, we linearize the general nonlinear system.
Indeed, linearization has become the norm in analyzing diverse networked systems \cite{Liu2011,Rajapakse2011pnas,Gao16Nature,Yan2017Nature,Gu2017Traj,Bassett18NP} due to several reasons.
One is that the empirical nonlinearity and the related parameters are hard to quantify and to estimate.
Another one is governed by a lemma that if the linearized system of a nonlinear dynamics is controllable along a specific trajectory, the nonlinear system is also controllable along the same trajectory \cite{LinearizationBook}.
The basic theoretical laws and insights of the practical control trajectory from initial to final state uncovered here facilitate the implementation of actual control in various empirical systems.
And it is worth further investigating for generalized scenarios of general nonlinear dynamics \cite{Khalil2002a} of static networks or temporal networks \cite{Holme2012,NaokiBook,Li2017,Li2017ConEng}.


\newpage

\end{document}